\documentclass[a4paper,11pt,leqno]{article}
\usepackage{amssymb}
\newcommand{\Tr}{\mathrm{Tr}}
\newcommand{\spde}{stochastic partial differential equation}
\newcommand{\ds}{\displaystyle}
\newcommand{\dd}{\mathrm{d}}

\newcommand{\der}[3]{\displaystyle{\frac{\partial^{#3} #1}{\partial #2^{#3}}}}

\newcommand{\R}{\mathbb R}
\newcommand{\E}{\mathbb E}
\newcommand{\N}{\mathbb N}

\newcommand{\ep}{\varepsilon}
\newcommand{\demi}{\frac{1}{2}}
\newcommand{\proba}[3]{(#1,{\cal #2},\mathbb{#3},\{{\cal #2}_t\}_{t \geq 0})}

\def\ds{\displaystyle}

\newtheorem{Theorem}{Theorem}[section]

\newtheorem{Proposition}[Theorem]{Proposition}
\newtheorem{Lemma}[Theorem]{Lemma}

\newtheorem{Example}[Theorem]{Example}

\def\qed{\hbox{\hskip 6pt\vrule width6pt height7pt
depth1pt  \hskip1pt}\bigskip}

\title{Weak order for the discretization of the stochastic heat equation.}
\author{
Arnaud Debussche\thanks{IRMAR et ENS de Cachan, antenne de Bretagne, Campus de Ker Lann, 
avenue Robert Schumann,
35170 BRUZ, FRANCE. {\tt arnaud.debussche@bretagne.ens-cachan.fr}} 
\and 
Jacques Printems\thanks{
Laboratoire d'Analyse et de Math\'ematiques Appliqu\'ees, CNRS UMR 8050, 
Universit\'e de Paris XII, 61, avenue du G\'en\'eral de Gaulle, 94010 CR\'ETEIL, FRANCE. {\tt printems@univ-paris12.fr}}}
\date{}

\begin{document}
\maketitle
\begin{abstract}
In this paper we study the approximation of the distribution of $X_t$ Hilbert--valued stochastic process solution of 
a linear parabolic stochastic partial differential equation written in an abstract form as 
$$
dX_t+AX_t \, dt = Q^{1/2} d W_t, \quad X_0=x \in H, \quad t\in[0,T],
$$  
driven by a Gaussian space time noise whose covariance operator $Q$ is given. We assume 
that $A^{-\alpha}$ is a finite trace operator for some $\alpha>0$ and 
that $Q$ is bounded from $H$ into $D(A^\beta)$ for some $\beta\geq 0$. It is not required 
to be nuclear or to commute with $A$.

The discretization is achieved thanks to finite element methods 
in space (parameter $h>0$) and implicit Euler schemes in time (parameter $\Delta t=T/N$). We define a discrete solution $X^n_h$ and 
for suitable functions $\varphi$ defined on $H$, we show that
$$
|\E \, \varphi(X^N_h) - \E \, \varphi(X_T) | =  O(h^{2\gamma} + \Delta t^\gamma)
$$
\noindent where $\gamma<1- \alpha + \beta$. Let us note that 
as in the finite dimensional case the rate of convergence is twice the one for pathwise 
approximations.
\end{abstract}



\section{Introduction}

In this article, we study the convergence of the distributions of numerical approximations 
of the solutions of the solutions of a large class of linear parabolic 
stochastic partial differential equations. 
The numerical analysis of stochastic partial differential equations has been recently the subject 
of many articles. (See among others \cite{allen-novosel-zhang}, \cite{davie-gaine}, 
\cite{greksch-kloeden}, \cite{gyongy1}, 
\cite{gyongy2}, \cite{gyongy-millet1}, \cite{gyongy-millet2},
\cite{gyongy-millet3}, \cite{gyongy-nualart}, \cite{hausenblas1},
\cite{hausenblas2}, \cite{lord-rougemont}, \cite{millet-morien}, \cite{printems}, \cite{shardlow}, 
\cite{walsh}, \cite{yan1}, 
\cite{yan2}). In all these papers, the aim is to give estimate on the strong order of convergence
for a numerical scheme. In other words, on the order of pathwise convergence. 
It is well known that in the case of stochastic differential equations in finite
dimension, the so-called weak order is much better than the strong order. The weak order is
the order of convergence of the law of approximations to the true solution. For instance, the standard
Euler scheme is of strong order $1/2$ for the approximation of a stochastic differential equation while
the weak order is $1$. A basic tool to study the weak order is the Kolmogorov equation associated to 
the stochastic equation (see \cite{kloeden-platten}, \cite{milstein}, \cite{milstein-tretyakov}
\cite{talay}).

In infinite dimension, this problem has been studied in much less articles. In 
\cite{buckwar-shardlow}, the case of a stochastic delay equations is studied.
To our knowledge, only \cite{dB-D}, \cite{haus3} consider stochastic partial differential equations. In 
\cite{dB-D} 
the nonlinear Schr\"odinger equations is considered. In the present article, we consider the 
case of the full discretization of a parabolic equation. We restrict our attention to 
a linear equation with additive noise which contains several difficulties. The general case 
of semilinear equations with state dependent noise present further difficulties and 
will be treated in a forthcoming article. This case is treated in \cite{haus3} but there only finite
dimensional functional of the solution are used and the finite dimensional method can be used.

Note that there are essential differences between the equations treated in \cite{buckwar-shardlow}
and \cite{dB-D}. Indeed, no spatial difference operator appear in a delay equation. In 
the case of the Schr\"odinger equation the linear evolution defines a group and it is possible 
to get rid of the differential operator by inverting the group. Furthermore, in 
 \cite{dB-D}, the data are assumed to be very smooth. 
 
In this article, we get rid of the differential operator by a similar trick as in \cite{dB-D}. However, 
since the linear evolution operator is not invertible, this introduces extra difficulties. Moreover, 
we consider a full discretization using an implicit Euler scheme and finite elements for the 
spatial discretization. We give estimate of the weak order of convergence with minimal 
regularity assumptions on the data. In fact we show that as in the finite dimensional case the weak
order is twice the strong order of convergence, both in time and space.

In dimension $d=1,2,3$, let us consider the following stochastic partial differential equation:
%$$
%\der{u(x,t)}{t}{} + (- \Delta)^{s/2} u(x,t) = \int_{\cal O} q(x-y) \dd W(y,t), 
%$$
%\begin{equation} \label{eq:intro0}
%\der{u(x,t)}{t}{} + (- \Delta)^{s/2} u(x,t) = \int_{\cal O} q(x-y) \xi(y,t) \dd y, 
%\end{equation}
\begin{equation} \label{eq:intro0}
\der{u(x,t)}{t}{} - \Delta u(x,t) = \dot{\eta}(x,t),
\end{equation}
\noindent where $x\in {\cal O}$, a bounded open set of $\R^d$, and $t\in]0,T]$, with Dirichlet boundary 
conditions and initial data and $\dot{\eta}=\der{\eta}{t}{}$ with $\eta$ denotes a real valued 
Gaussian process. It is convenient to use an abstract framework to describe the noise more precisely.
Let $W$ be a cylindrical Wiener process on $L^2({\cal O})$, in other words $\ds{\frac{\partial W}{\partial t}}$
is the space time white noise. Equivalently, given an orthonormal basis of $L^2({\cal O})$, 
$W$ has the following expansion
$$
W(t) = \sum_{i\in \N} \beta_i(t) e_i
$$
where $(\beta_i)_{i\in•À?\N}$ is a family of independent standard brownian motions 
(See section 2.4 below). We consider noises of the form 
$
\eta(t) = Q^{1/2} W(t)
$
where $Q$ is a non negative symmetric bounded linear operator on $L^2({\cal O})$.
For  example{\footnote{The equations describing this example are formal, it is not difficult to give
a rigourous meaning. This is not important in our context. Note that $q$ does not need to be 
a function and a distribution is allowed.}, given a function $q$ defined on ${\cal O}$, we can take 
$$
\eta(x,t) = \int_{\cal O} q(x-y) W(y,t) \dd y,
$$
Then the process $\eta$ has the following correlation function:
$$
\E \, \eta(x,t) \eta(y,s) = c(x-y) (t\wedge s) \quad \mbox{with   } 
c(r) = \int_{\cal O} q(z+r) q(z) \, \dd z.
$$
The operator $Q$ is then given by  
$
Q f(x) = \int_{\cal O} c(x-y) f(y) \, \dd y.
$
Note that if the $q$ is the Dirac mass at $0$, $\eta=W$ and $Q=I$.

Let us also set $A=-\Delta$, $D(A)=H^2({\cal O})\cap H^1_0({\cal O})$ 
and $H=L^2({\cal O})$. Then $A:D(A)\rightarrow H$ can be seen as an 
unbounded operator on $H$ with domain  $D(A)$.
Our main assumption concerning $Q$ is that $A^{\sigma/2}Q$ defines a bounded operator on 
$L^2({\cal O})$ with  
$\sigma>-1/2$ if $d=1$, $\sigma>0$ if $d=2$ and $\sigma>1/2$ if $d=3$. 
In the example above this amounts to require that 
$(-\Delta)^{\sigma/2} q \in L^2({\cal O})$.
It is well known that theses conditions are sufficient to ensure the existence of continuous 
solutions of (\ref{eq:intro0}). 

If we write $u(t) = u(\cdot,t)$ seen as a 
$H$--valued stochastic process then (\ref{eq:intro0}) can be rewritten under the abstract Ito 
form
\begin{equation} \label{eq:abs}
\dd u(t) + A u(t) \, \dd t = Q^{1/2} \, \dd W(t).
\end{equation}
In this article, we consider such an abstract equation and study the approximation of the law of the solutions of (\ref{eq:abs}) by means of finite elements of the 
distribution of $u$ in $H$. Let $\{V_h\}_{h\geq 0}$ be a family of finite dimensional subspaces 
of $D(A^{1/2})$. Let $N\geq 1$ an integer and $\Delta t = T/N$. The numerical scheme is given by 
\begin{equation}\label{eq:scheme}
(u^{n+1}_h - u^n_h,v_h) + \Delta t(A u^{n+\theta}_h,v_h)  = \sqrt{\Delta t} \, (Q^{1/2} \chi^{n+1},v_h),
\end{equation}
\noindent for any $v_h \in V_h$, where $\sqrt{\Delta t}\,  \chi^{n+1} = W((n+1)\Delta t) - W(n\Delta t)$ is the noise 
increment and where $(\cdot,\cdot)$ is the inner product of $H$. The unknown is approximated at time $n \Delta t$, $0\leq n \leq N$ by $u_h^n \in V_h$. 
In (\ref{eq:scheme}), we have used the notation $u^{n+\theta}_h$ for $\theta u_h^{n+1} + 
(1-\theta)u_h^n$ for some $\theta \in [0,1]$. We prove that an error estimate of the following form
$$
|\E(\varphi(u(n\Delta t))-\E(\varphi(u_h^n))|\le c (h^{2\gamma} + \Delta t^{\gamma})
$$
for any function $\varphi$ which is $C^2$ and bounded on $H$. With  the above notation,
$\gamma$ is required to be strictly less than $1-d/2+\sigma/2$. This is 
exactly twice the strong order (see \cite{yan1}, \cite{walsh}). If $d=1$ and $\sigma=0$,
the condition is $\gamma<1/2$ and we obtain a weak order which $1/2$ in time and $1$ in 
space.

%Concerning the deterministic case ($q=0$), the results about the rate of convergence are 
%classical. Nevertheless, some previous works (see \cite{crouzeix_et_al_89,leroux}) show that the infinite 
%dimensional setting implies 
%$$
%\theta >1/2
%$$
%\noindent for the convergence.

%Note that since the unknowns are $H$--valued random variables, various 
%type of convergence can be considered, either we consider pathwise approximation (strong order) 
%or approximation 
%of the distribution (weak order). 
%In the finite dimensional case, discretization in time of 
%stochastic 
%differential equations is now well understood . Typically, 
%under suitable conditions on the coefficients, the Euler scheme is of strong order 1/2 and 
%of weak order 1.
%Rsults about the rate of convergence in the infinite dimensional case is not developped as in 
%the finite dimensional case (see \cite{kloeden,milstein,talay})

%The main result of this paper is to estimate the error between $\E \varphi(u(T))$ and 
%$\E \varphi(u_h^N))$ for some functions $\varphi:H\rightarrow \R$ in function of $\Delta t$ 
%and $h$. Note that contrary to the 

\section{Preliminaries} \label{sec:notations}
\subsection{Functional spaces.}
It is convenient to change the notations and rewrite the unknown of (\ref{eq:abs}) 
as $X$. We thus consider the following stochastic partial differential equation written 
in the  abstract form
\begin{equation} \label{eq:0}
\dd X_t + A X_t \, \dd t = Q^{1/2} \dd W_t, \quad X_0 = x \in H, \quad 0<t\leq T, 
\end{equation}
\noindent where $H$ is a Hilbert space whose inner product is denoted 
by $(\cdot,\cdot)$ and its associated norm by $|\cdot|$, the process $\{X_t\}_{t\in[0,T]}$ is 
an $H$--valued stochastic process, $A$ 
a non negative self-adjoint unbounded operator on $H$ whose domain $D(A)$ is dense in $H$ and 
compactly embedded in $H$, $Q$ a non negative symmetric operator on $H$ and $\{W_t\}_{t\in[0,T]}$ a 
cylindrical Wiener process on $H$ adapted to a given normal filtration $\{{\cal F}_t\}_{t\in [0,T]}$ in a given probability space 
$(\Omega,{\cal F},{\mathbb P})$.

It is  well known that there exists a sequence of nondecreasing positive real numbers $\{\lambda_n\}_{n\geq 1}$ 
together with $\{e_n\}_{n\geq 1}$ a Hilbertian basis of $H$ such that
$$
A e_n = \lambda_n e_n \quad \mbox{with   }\lim_{n\rightarrow +\infty}\lambda_n =  + \infty.
$$
We set for any $s\geq 0$, 
$$
D(A^{s/2}) = \left \{ u =\sum_{n=1}^{+\infty} u_n e_n  \in H \quad \mbox{such that   } 
\sum_{n=1}^{+\infty}\lambda_n^{s} u_n^2 < +\infty \right \},
$$
\noindent and 
$$
A^s u = \sum_{n\geq 1} \lambda^s_n u_n e_n, \quad \forall u \in D(A^s).
$$
It is clear that $D(A^{s/2})$ endowed with the norm $u\mapsto \|u\|_s :=  |A^{s/2} u| $ 
is a Hilbert space. We define also $D(A^{-s/2})$ with $s\geq 0$ as the completed space 
of $H$ for the topology induced by 
the norm $\|u\|_{-s}=\sum_{n\geq 1} \lambda^{-s}_n u_n^2$. In this case $D(A^{-s/2})$ can be identified with the topological dual of $D(A^{s/2})$, i.e. the space of the 
linear forms on $D(A^{s/2})$ which are continuous with respect to the topology induced by the norm  $\|\cdot\|_s$.

Moreover, theses spaces can be obtained by interpolation between them. Indeed, for any reals $s_1 \leq  s \leq  s_2$, one has the continuous embeddings $D(A^{s_1/2}) \subset D(A^{s/2}) \subset D(A^{s_2/2})$ and by H\"older inequality
\begin{equation} \label{eq:DAs_interp}
\| u \|_s \leq  \, \| u \|_{s_1}^{1-\lambda} \|u \|_{s_2}^{\lambda}, \quad s = (1-\lambda) s_1 + \lambda s_2,  
\end{equation}
\noindent for any $u \in D(A^{s_2/2})$.

We denote by $\|\cdot\|_X$ the norm of a Banach space $X$. If $X$ and $Y$ 
denote two Banach spaces, we
denote by ${\cal L}(X,Y)$ the Banach space of 
bounded linear operators from $X$ into $Y$ endowed with the norm 
$\|B\|_{{\cal L}(X,Y)}=\sup_{x\in X} \|Bx\|_Y/\|x\|_X$. When $X=Y$, we use the shorter notation 
${\cal L}(X)$. 

If $L \in {\cal L}(H)$ is a nuclear operator, $\Tr(L)$ denotes the trace of the operator $L$, i.e. 
$$
\Tr(L) = \sum_{i\geq 1} (L e_i,e_i) \leq + \infty.
$$
It is well known that the previous definition does not depend on the choice of the Hilbertian 
basis. Moreover, the following properties hold
\begin{equation} \label{eq:trace_prop1}
\Tr{}(LM) = \Tr{}(ML), \quad \mbox{for any  } L,M \in {\cal L}(H),
\end{equation}
\noindent and 
\begin{equation} \label{eq:trace_prop2}
\Tr{}(LM) \leq \Tr{}(L) \|M\|_{{\cal L}(H)}, \quad \mbox{for any  } L \in {\cal L}_+(H), \;M \in {\cal L}(H),
\end{equation}
\noindent where ${\cal L}_+(H)$ denotes the set  non negative bounded linear operators on $H$.

Hilbert-Schmidt operators play also an important role. Given two Hilbert spaces $K_1,\; K_2$, 
an operator $L \in {\cal L}(K_1,K_2)$ is Hilbert-Schmidt if $L^*L$ is a  nuclear operator on $K_1$ or 
equivalently if $LL^*$ is nuclear on $K_2$. We denote by ${\cal L}_2(K_1,K_2)$ the space of
such operators. It is a Hilbert space for the norm
$$
\|L\|_{{\cal L}_2(K_1,K_2)}= \left( \Tr{} (L^*L)\right)^{1/2}=  \left( \Tr{} (L L^*)\right)^{1/2}.
$$
It is classical that, given four Hilbert spaces $K_1,\; K_2,\; K_3,\; K_4$,  if $L\in {\cal L}_2(K_2,K_3)$,
$M\in {\cal•À?L}(K_1,K_2)$, $N\in {\cal L}(K_3,K_4)$ then $NLM\in {\cal L}_2(K_1,K_4)$ and 
\begin{equation} \label{eq:HS}
\| NLM\|_{{\cal L}_2(K_1,K_4)}\le \| N\|_{{\cal L}(K_3,K_4)} 
\|L\|_{{\cal L}_2(K_2,K_3)}\|M\|_{{\cal•À?L}(K_1,K_2)}.
\end{equation}
See \cite{dpz}, appendix C, or \cite{gohberg-krejn} for more details on nuclear and Hilbert-Schmidt 
operators.

If $X$ is a Banach space, we denote by ${\cal C}_b(H;X)$ the Banach space of $X$-valued, continuous and bounded functions on $H$. We also denote by $C^k_b(H)$ the space of $k$-times continuously
differentiable real valued functions on $H$. The first order differential of a function 
$\varphi \in C^1_b(H)$ is identified with its gradient and is then considered as an 
element of $C_b(H;H)$. It is denoted
by $D\varphi$. Similarly, the second order differential of a function 
$\varphi \in C^2_b(H)$ is seen as  a function from $H$ into the Banach space ${\cal L}(H)$
and is denoted by $D^2\varphi$.

%Let us describe in what context such spaces can occur. If $v : H\rightarrow \R$ denotes a real valued function defined on the Hilbert space $H$, we  
%denote, when it exists, by $Dv(a)$ its Fr\'echet-differential at $a\in H$, i.e. the element 
%of $H$ such that 
%$$
%\lim_{|k|\rightarrow 0} \frac{|v(a+k) - v(a) - (Dv(a),k)|}{|k|} = 0.
%$$
%When it exists, $D^2v(a)$ denotes the Fr\'echet-differential of $a\mapsto Dv(a)$ at $a\in H$. It 
%will be understood  as an element of ${\cal L}(H)$. Hence, the function $D^2v$ can be seen as an application 
%from $H$ into the Banach space ${\cal L}(H)$. 

\subsection{The deterministic stationary problem}

We need some classical results on the deterministic stationary version of 
(\ref{eq:0}). In this case, special attention has to be paid to the space $V=D(A^{1/2})\subset H$. It is a Hilbert space whose embedding into $H$ is dense and continuous. Its inner product is denoted by $((\cdot,\cdot))$. We have
$$
\begin{array}{rcll}
((u,v)) & = & (A^{1/2}u, A^{1/2}v), & \mbox{for any   }u\in V,  \; v \in V \\
            & = & (Au,v),                         & \mbox{for any   }u \in D(A), \; v \in H.
%        & = & \sum_{i\geq 1}^{+\infty} \lambda_i u_i v_i
\end{array}
$$
Then by a density argument and the uniqueness of the Riesz representation (in $V$) we conclude that $A$ is invertible 
from $V$ into $V'=D(A^{-1/2})$ or from $D(A)$ into $H$. We will set $T=A^{-1}$ its inverse. It is 
bounded and positive on $H$ and on $V$.

For any $f \in H$, $u=T f$ is by definition the unique solution of the following problem
\begin{equation} \label{eq:defT}
u\in V, \quad ((u,v)) = (f,v), \quad \mbox{for any   } v \in V.
\end{equation}

Let $\{V_h\}_{h>0}$ be a family of finite dimensional subspaces of $V$ parametrized by a small 
parameter $h>0$. For any $h>0$, we denote by 
$P_h$ (resp. $\Pi_h$) the orthogonal projector from $H$ (resp. $V$) onto $V_h$ with respect to 
the inner product $(\cdot,\cdot)$ (resp. $((\cdot,\cdot))$).

For any $h>0$, we denote by $A_h$ the 
linear bounded operator from $V_h$ into $V_h$ defined by 
\begin{equation} \label{eq:defAh}
((u_h,v_h)) = (A_h u_h,v_h) = (A u_h,v_h) \quad \mbox{for any  } u_h \in V_h, \; v_h \in V_h.
\end{equation}
It is clear that $A_h:V_h \rightarrow V_h$ is also invertible. Its inverse is denoted by $T_h$. 
For any $f\in H$, $u_h = T_h f$ is by definition the solution of the following problem :
\begin{equation} \label{eq:defTh}
u_h \in V_h, \quad ((u_h,v_h)) = (f,v_h) = (P_h f,v_h),\quad \mbox{for any  }v_h \in V_h.
\end{equation}

It is also clear that $A_h$ and $T_h$ are positive definite symmetric bounded linear operators on $V_h$. We 
denote by $\{\lambda_{i,h}\}_{1\leq i \leq I(h)}$ the sequence of its nonincreasing positive 
eigenvalues and $\{e_{i,h}\}_{1\leq i \leq I(h)}$ 
the associated orthonormal basis of $V_h$ of its eigenvectors. Again, by H\"older inequality,
 $A_h$ satisfies the following interpolation inequality
 \begin{equation}
 \label{interp_h}
 |A_h^su_h|\le |A_h^{s_1}u_h|^{\lambda}|A_h^{s_2}u_h|^{1-\lambda},\; u_h\in V_h,\;
 s=\lambda s_1 +(1-\lambda)s_2.
 \end{equation}

The consequences of (\ref{eq:defAh}) are summarized in the following Lemma. 
\begin{Lemma}
Let $A_h \in {\cal L}(V_h)$ defined in (\ref{eq:defAh}). Let $T$ and $T_h$ defined in (\ref{eq:defT}) and (\ref{eq:defTh}). Then
%\begin{equation} \label{eq:lambda}
%\lambda_{i,h} \geq \lambda_i, \qquad 1\leq i \leq I(h). \; UTILE ?
%\end{equation}
%Moreover 
the following hold for any $w_h \in V_h$ and $v\in V$:
\begin{eqnarray} \label{eq:cea}
T_h P_h & = & \Pi_h T, \\
\label{eq:interp_A}
|A^{1/2} w_h|    & = & |A_h^{1/2} w_h| \\
\label{eq:interp_T-1}
|T_h^{1/2} w_h| & = & |T^{1/2} w_h|,\\
\label{eq:interp_T-2}
|T_h^{1/2} P_h v| & \leq & |T^{1/2} v |.
\end{eqnarray}
\end{Lemma}
\noindent {\bf Proof}

%Using the min-max principle of Courant--Fischer,  (\ref{eq:lambda}) is easily shown. Indeed, let $i\geq 1$ and
% $W_{i,h}$ a $i$-dimensional subspace of $V_h$. Then, (\ref{eq:defAh}) implies that for any 
%$v_h \in W_{i,h}$, one has 
%$$
%\max_{v_h \in W_{i,h},\; v_h\neq 0} \frac{(A_h v_h,v_h)}{|v_h|^2} = \max_{v_h \in W_{i,h},\; v_h\neq 0} \frac{(A v_h,v_h)}{|v_h|^2}.
%$$
%Now, since $V_h\subset V$, $W_{i,h}$ is itself a $i$-dimensional subspace of $V$, the right
%hand side of the previous equality is 
%greater or equal than $\lambda_i$. We conclude by taking  the minimum over all the $i$-dimensional subspace of $V_h$.

Now let $f\in H$. We consider the two solutions $u$ and $u_h$ of (\ref{eq:defT}) and (\ref{eq:defTh}). Since $V_h \subset V$, we can 
write (\ref{eq:defT}) with $v_h \in V_h$. Then, substracting we get $((u-u_h,v_h))=0$. Hence, $u_h=\Pi_h u$ the $V$-orthogonal 
projection of $u$ onto $V_h$, i.e. $T_h P_h f = \Pi_h T f$.

Equation (\ref{eq:interp_A}) follows immediately from the definition (\ref{eq:defAh}) of $A_h$. We now prove (\ref{eq:interp_T-2}). Equations (\ref{eq:cea}) and (\ref{eq:interp_A}) imply $|A_h^{-1/2}P_h v| = \| T_h P_h v\| = \| \Pi_h Tv\| \leq  
\|T v\| =|A^{-1/2} v|$, since $\Pi_h:V\rightarrow V_h$ is an orthogonal projection for the inner product 
$\|\cdot\|$.

As regards (\ref{eq:interp_T-1}), on one hand (\ref{eq:interp_T-2}) with $v=w_h \in V_h \subset V$ gives the first inequalty 
$|T_h^{1/2} w_h| \leq |T^{1/2} w_h|$. On the other hand, by (\ref{eq:defAh}),
$$
|(A_h u_h,v_h)|  = |(A u_h,v_h)| \leq |A^{1/2} u_h| |A^{1/2} v_h|= |A^{1/2} u_h| | A^{1/2} v_h|. 
$$
So $A_h u_h$ can be considered as a continuous linear form on $D(A^{1/2})$, i.e. belongs to $D(A^{-1/2})$, and
$$
| A^{-1/2} A_h u_h | \leq |A^{1/2}u_h| = |A_h^{1/2} u_h|.
$$
\noindent Taking $u_h = A_h^{-1} w_h$ gives $|A^{-1/2} w_h| \leq |A_h^{-1/2} w_h|$. Eq. (\ref{eq:interp_T-1}) follows.
 \hfill \qed

Our main assumptions concerning the spaces $V_h$ is that the corresponding linear elliptic problem (\ref{eq:defTh}) admits an 
$O(h^r)$ error estimates in $H$ and $O(h^{r-1})$ in $V$ for some $r\geq 2$. It is classical to verify that these estimates hold if we suppose 
that $\Pi_h$ satisfies 
for some constant $\kappa_0>0$,
\begin{eqnarray}\label{eq:FEM_H}
|\Pi_h v - v|                    & \leq & \kappa_0 \, h^s |A^{s/2} v|, \quad 1\leq s \leq r,\\
|A^{1/2} (\Pi_h w - w)| & \leq & \kappa_0 \, h^{s'-1} |A^{s'/2} w|,\quad 1 \leq s' \leq r-1,\label{eq:FEM_V}
\end{eqnarray}
\noindent where $v\in D(A^{s/2})$ and $w\in D(A^{s'/2})$.

Finite elements satisfying these conditions are for example $P_k$ triangular elements on a 
polygonal domain or $Q_k$ rectangular finite element on a rectangular domain provided
$k\ge 1$. Approximation by splines can also be considered. (See \cite{ciarlet78}, \cite{strang-fix}).

\subsection{The deterministic evolution problem}

We recall now some results about the spatial discretization of the solution of the deterministic linear 
parabolic evolution equation:
\begin{equation} \label{eq:heat}
\der{u(t)}{t}{} + Au(t) = 0, \quad u(0) = y,
\end{equation}
\noindent by the finite dimensional one
$$
\der{u_h(t)}{t}{} + A_h u_h(t) = 0, \quad u_h(0) = P_h y \in V_h.
$$

It is well known that, under our assumptions, (\ref{eq:heat}) defines a contraction semi-group on $H$ denoted by $S(t)=e^{-tA}$ for any $t\geq 0$. 
Its solution 
can be read as $u(t) = S(t) y$ where $t\geq 0$. The main properties of $S(t)$ (contraction, regularization) are summed up below:
\begin{equation} \label{eq:contraction}
|e^{-tA} x| \leq |x|, \quad \mbox{for any   } x \in H,
\end{equation}
\noindent and
\begin{equation}\label{eq:regular0}
| A^s e^{-t A} x| \leq C(s) t^{-s} |x|,
\end{equation}
\noindent for any $t>0$, $s\geq 0$ and $x\in H$. Such a property is based on the definition of $A^s$ and the following well known inequality
\begin{equation}\label{eq:regular}
\sup_{x\geq 0} x^\ep e^{-tx} \leq C(\ep) t^{-\ep}, \quad \mbox{for any  } t > 0.
\end{equation}

In the same manner, we denote by $S_h(t)$ or $e^{-t A_h}$ 
the semi-group on $V_h$ such that 
$u_h(t) = S_h(t) P_h y$, for any $t\geq 0$. 

We have various types of convergence of $u_h$ towards $u$ depending on the regularity of the initial data $y$. The optimal rates of convergences remain the same as in the corresponding stationary problems (see (\ref{eq:FEM_H})--(\ref{eq:FEM_V})).  The estimates are not uniform in time near $t=0$ 
since the regularization of $S(t)$ is used to prove them. 
The following Lemma gives two classical properties needed in this article.

\begin{Lemma} \label{lem:smooth_or_not}
Let $r\geq 2$ be such that (\ref{eq:FEM_H}) and (\ref{eq:FEM_V}) hold and $q$, $q'$, $s$, $s'$ such that $0\leq s \leq q \leq r$, 
$s'\geq 0$, $1\leq q' +s' \leq r-1$ and $q'<2$. 
Then there exists constants $\kappa_i>0$, $i=1,2$ independent on $h$ such that for any time $t>0$, one has: 
\begin{eqnarray} 
\label{eq:smooth}
\|S_h(t)P_h - S(t) \|_{{\cal L}(D(A^{s/2}),H)} & \leq  & \kappa_1 \, h^q \, t^{-(q-s)/2}, \\
\label{eq:non-smooth0}
\|S_h(t)P_h - S(t)\|_{{\cal L}(D(A^{s'/2}),D(A^{1/2}))}  & \leq  & \kappa_2 \, h^{q'+s'} \, t^{-(q'+1)/2}.
%\label{eq:non-smooth}
%\|S_h(t)P_h - S(t)\|_{{\cal L}(H,D(A^{\delta/2}))} & \leq & \kappa_3 \, h^{q-\delta} t^{-q/2},
\end{eqnarray}
%\noindent for $\delta \in [0,1]$.
\end{Lemma}

The proof of (\ref{eq:smooth}) can be found in \cite{bramble_et_al_77} (see also \cite{thomee97}, Theorem 3.5, p. 45). The proof of (\ref{eq:non-smooth0}) can be found in \cite{johnson87}, Theorem 4.1, p. 342 (with $f=0$)). In fact, we use only $s=s'=0$ and $q'=1$ 
below.

%Finally, (\ref{eq:non-smooth}) 
%follows from the interpolation between (\ref{eq:smooth}) with $s=0$ and (\ref{eq:non-smooth0}) with $s'=0$ and $q'=q-1$.

%IL FAUT DONC SUPPOSER $q\ge 2$.

\subsection{Infinite dimensional stochastic integrals}

In this section, we recall basic results on the stochastic integral with respect to the 
cylindrical Wiener process $W_t$. More details can be found for instance in \cite{dpz}.

It is well known that  $W_t$ has the following expansion 
$$
W_t = \sum_{i=1}^{+\infty} \beta_i(t) e_i,
$$
\noindent where $\{\beta_i\}_{i\geq 1}$ denotes a family of real valued mutually independent Brownian motions on 
$\proba{\Omega}{F}{P}$. The sum does not converge in $H$ and this reflects the bad regularity 
property of the cylindrical Wiener process. However, it converges a.s. and in $L^p(\Omega;U)$, 
$p\ge 1$,
for any space $U$ such that $H\subset U$ with a Hilbert-Schmidt embedding. If $H=L^2({\cal O})$,
${\cal O}\subset \R^d$ open and bouded, we can take $U=H^{-s}({\cal O})$, $s>d/2$.

Such a Wiener process can be characterized by
$$
\E \, (W_t, u) (W_s,v) = \min(t,s) (u,v)
$$
\noindent for any $t,s\geq 0$ and $u,v \in H$.

Given any predictable operator valued function $t\mapsto \Phi(t),\; t\in [0,T]$, it is possible to define 
$\int_0^T \Phi(s)dW(s)$ in a Hilbert space $K$ if $\Phi$ takes values in ${\cal L}_2(H,K)$ 
and 
$\int_0^T \|\Phi(s)\|_{{\cal L}_2(H,K)}^2 ds <\infty$ a.s. In this case $\int_0^T \Phi(s)dW(s)$ is 
a well defined random variable with values in $K$ and 
$$
\int_0^T \Phi(s)dW(s)= \sum_{i=1}^\infty \int_0^T \Phi(s)e_i d\beta_i(s).
$$
Moreover, if 
$\E\left(\int_0^T \|\Phi(s)\|_{{\cal L}_2(H,K)}^2 ds\right) <\infty$, then 
$$
\E\left(\int_0^T \Phi(s)dW(s)\right) =0,
$$
and 
$$
\E\left(\left(\int_0^T \Phi(s)dW(s)\right)^2\right)= \E\left(\int_0^T \|\Phi(s)\|_{{\cal L}_2(H,K)}^2 ds\right).
$$
We will consider below expressions of the form $\int_0^t \psi(s) Q^{1/2}  dW(s)$. These are then
square integrable random variables in  $H$ with zero average if 
$$
\E\left(\int_0^T \|\psi(s)Q^{1/2}\|_{{\cal L}_2(H,K)}^2 ds\right)=
\E\int_0^T \Tr{}\left(\psi^*(s) Q\psi(s)\right) ds<\infty.
$$
The solution of equation (\ref{eq:0}) can be written explicitly in terms of stochastic integrals. In order
that these are well defined, we assume throughout this paper that there exists real numbers 
$\alpha >0$ 
and $\min(\alpha-1,0)\leq \beta\leq \alpha$ such that
\begin{equation} \label{eq:trace}
\sum_{n=1}^\infty \lambda_n^{-\alpha} = \|A^{-\alpha/2}\|_{{\cal L}_2(H)}^2=\Tr(A^{-\alpha}) < +\infty.
\end{equation}
\noindent and
\begin{equation} \label{eq:Q}
A^\beta Q \in {\cal L}(H).
\end{equation}
Condition (\ref{eq:Q}) implies that $Q$ is a bounded operator from $H$ into $D(A^\beta)$.
By interpolation, we deduce immediately that for any $\lambda\in [0,1]$, $A^{\lambda\beta} Q^\lambda \in {\cal L}(H)$ and 
\begin{equation} \label{eq:Q_lambda}
\|A^{\lambda\beta} Q^\lambda \|_{\cal L}(H)\le \|A^\beta Q\|_{{\cal L}(H)}^\lambda.
\end{equation}
\begin{Example}
If one considers the equations described in the introduction where $A$ is the Laplace 
operator with Dirichlet boundary conditions,  it is well known that (\ref{eq:trace}) holds 
for $\alpha>d/2$.
\end{Example}

We have the following result.
\begin{Proposition}\label{prop:X}
Assume that (\ref{eq:trace}), (\ref{eq:Q}) hold and
\begin{equation}\label{eq:order}
1 - \alpha + \beta >0.
\end{equation}
Then there exists a unique Gaussian stochastic process which is the weak solution (in the PDE sense) of (\ref{eq:0}) continuous in time with values in $L^2(\Omega,H)$. It is given by the formula which holds a.s. in $H$:
$$
X_t = e^{-tA} x + \int_0^t e^{-(t-s)A} Q^{1/2} \, \dd W_s = e^{-tA} x + 
\sum_{i=1}^{+\infty}\left (\int_0^t e^{-(t-s)\lambda_i} \dd \beta_i(s) \right ) Q^{1/2} e_i.
$$
\end{Proposition}
\noindent {\bf Proof}

By  Theorem 5.4 p. 121 in \cite{dpz}, 
it is sufficient to see that  the stochastic integral make sense in $H$, i.e. 
$$
\int_0^t \|e^{-(t-s)A} Q^{1/2}\|_{{\cal L}_2(H)}^2ds = \int_0^t \Tr{} \left (  e^{-(t-s)A} Q e^{-(t-s)A} \right ) ds <\infty,
$$
\noindent for any $t \in [0,T]$. 
We use (\ref{eq:HS}) to estimate the Hilbert-Schmidt norm:
$$
\begin{array}{ll}
\|e^{-(t-s)A} Q^{1/2}\|_{{\cal L}_2(H)} &\le \| A^{\beta/2}Q^{1/2}\|_{{\cal L}(H)} 
\|A^{-\alpha/2}\|_{{\cal L}_2(H)}
\|e^{-(t-s)A} A^{(-\beta+\alpha)/2}\|_{{\cal L}(H)}\\
&\le c (t-s)^{-1/2(-\beta+\alpha)}
\end{array}
$$
by (\ref{eq:regular0}), (\ref{eq:trace}) and (\ref{eq:Q_lambda}). The conclusion follows since 
$-\beta+\alpha <1$.
\qed

%We choose $0<\ep<1-\alpha+\beta$ and write for $s\in[0,t[$:
%\begin{eqnarray*}
%\Tr{} \left (  e^{-(t-s)A} Q e^{-(t-s)A} \right ) & = & \sum_{i=1}^{+\infty} (e^{-2(t-s)A} Qe_i,e_i)  \\
%   & = & \sum_{i=1}^{+\infty} (A^{1-\ep} e^{-2(t-s)A} A^{\ep-1} Qe_i,e_i) \\
%   & \leq & \|A^{1-\ep} e^{-2(t-s)A} \|_{{\cal L}(H)} \sum_{i=1}^{+\infty} | (A^{-\beta+\ep-1} A^\beta Q e_i,e_i)| \\
%   & \leq & C(t-s)^{-(1-\ep)} \|A^\beta Q\|_{{\cal L}(H)} \sum_{i=1}^{+\infty} \lambda_i^{-\beta+\ep-1} \\
%   & \leq & C'(t-s)^{-(1-\ep)},
%\end{eqnarray*}
%\noindent where (\ref{eq:regular0}), (\ref{eq:trace}) and (\ref{eq:Q})  have been used.

%In the same way we prove easily by interpolation that 
%$Q \in {\cal L}(D(A^{-\beta+\lambda}),D(A^\lambda))$ for 
%any {\em real} $\lambda\in [0,\beta]$ and
%\begin{equation}
%\label{lem:Q}
%\|A^\lambda Q A^{-\beta+\lambda}\|_{{\cal L}(H)}=\|Q\|_{{\cal L}(D(A^{-\beta+\lambda}),D(A^\lambda))} \leq \|A^\beta Q\|_{{\cal L}(H)}.
%\end{equation}
%Indeed the result is clear for $\lambda =0$ and $\lambda =-\beta$.

%\begin{Lemma}\label{lem:Q}
%Let $Q$ be a symmetric operator on $H$ such that (\ref{eq:Q}) holds. Then $Q \in {\cal L}(D(A^{-\beta+\lambda}),D(A^\lambda))$ for 
%any {\em real} $\lambda\in [0,\beta]$ and
%$$
%\|Q\|_{{\cal L}(D(A^{-\beta+\lambda}),D(A^\lambda))} \leq \|A^\beta Q\|_{{\cal L}(H)}.
%$$
%\end{Lemma}
%The proof is very easy. Indeed, for $\lambda =0$ and $\lambda =-\beta$, the result is clear. The 
%general case follows by interpolation. 

%----------------------------------
%    MAIN RESULT AND PROOF
%----------------------------------
\section{Weak convergence of an implicit scheme.}
\subsection{Setting of the problem and main result.}
In this section, we state the weak approximation result 
on the full discretization of (\ref{eq:0}).
%of the following linear \spde{} 
%written in an abstract framework and in the It\^o sense
%\begin{equation} \label{eq:sto_heat}
%\dd X_t + A X_t = Q^{1/2} \dd W_t, \quad t\in[0,T],
%\end{equation}
%\noindent with the intial condition
%\begin{equation} \label{eq:ci}
%X_0 = x \in H,
%\end{equation}
%\noindent where the Hilbert $H$ and the linear operator $A$ and $Q$ have been introduced in Section 
%\ref{sec:notations}, 
%and where $\{W_t\}_{t\in[0,T]}$ denotes a cylindrical Wiener process on $H$. 

We first describe the numerical scheme. 
Let $N\geq 1$ be an integer and $\{V_h\}_{h>0}$ the family of finite 
element spaces introduced in Section \ref{sec:notations}. Let $\Delta t = T/N$ and 
$t_n = n\, \Delta t$, $0\leq n\leq N$. For any $h>0$ and any integer $n\leq N$, we seek for 
$X_h^n$, an approximation of $X_{t_n}$, such that for any $v_h$ in $V_h$:
\begin{equation} \label{eq:var_sto_heat}
(X_h^{n+1} - X_h^n,v_h) + \Delta t \, (A X_h^{n+\theta},v_h) = (Q^{1/2} W_{t_{n+1}} - Q^{1/2} W_{t_n},v_h), 
\end{equation}
\noindent with the initial condition
\begin{equation} \label{eq:var_ci}
(X_h^0,v_h) = (x,v_h), \quad \forall v_h \in V_h,
\end{equation}
\noindent where 
$$
X_h^{n+\theta} = \theta X_h^{n+1} + (1 - \theta) X_h^n,
$$
\noindent with
\begin{equation} \label{eq:implicit}
1/2 < \theta \leq 1.
\end{equation}
Recall that for $\theta \le 1/2$, the scheme is in general unstable and a CFL condition is necessary. 

Then (\ref{eq:var_sto_heat})--(\ref{eq:var_ci}) can be rewritten as
\begin{equation} \label{eq:disc_sto_heat}
X_h^{n+1} - X_h^n + \Delta t \, A_h X_h^{n+\theta}  = \sqrt{\Delta t} P_h Q^{1/2} \chi^{n+1}, 
\end{equation}
\begin{equation} \label{eq:disc_ci}
X_h^0 = P_h x,
\end{equation}
\noindent where 
$$
\chi^{n+1} = \frac{1}{\sqrt{\Delta t}} \left (W_{(n+1)\Delta t} - W_{n\Delta t} \right ),
$$
\noindent and where we recall that $P_h : H \rightarrow V_h$ is the $H$-orthogonal projector.  Hence 
$\{\chi^n\}_{n\geq 0}$ is a sequence of independent and identically distributed gaussian random variables. 
%---------------------
% THEOREME
%----------------------
The main result of this paper is stated below.
\begin{Theorem} \label{theo:weak}
Let $\varphi\in C^2_b(H)$, i.e.  a twice differentiable 
real valued functional defined on $H$ whose first and second derivatives are bounded. Let 
$\alpha>0$ and $\beta\geq 0$ be such that (\ref{eq:trace}), (\ref{eq:Q}) and  (\ref{eq:order}) hold. 
Let $T\geq 1$ and 
$\{X_t\}_{t\in [0,T]}$ be the $H$-valued 
stochastic process solution of (\ref{eq:0}) given by Proposition \ref{prop:X}. For any $N\geq 1$, let 
$\{X_h^n\}_{0\leq n\leq N}$ be the solution 
of the 
scheme (\ref{eq:disc_sto_heat})--(\ref{eq:disc_ci}). Then there 
exists a constant $C=C(T,\varphi)>0$ which does 
not depend on $h$ and $N$ such that 
for any $\gamma<1-\alpha+\beta\leq 1$, the following inequality holds
\begin{equation} \label{eq:weak}
\left | \E \, \varphi (X_h^N) - \E \, \varphi(X_T) \right | \leq C \left ( h^{2\gamma} + \Delta t^\gamma 
\right ),
\end{equation}
\noindent where $\Delta t = T/N \leq 1$. 
\end{Theorem}

%
% Preuve du Th•À?r•À?e 3.1
%
\subsection{Proof of Theorem \ref{theo:weak}.}

The scheme (\ref{eq:disc_sto_heat})--(\ref{eq:disc_ci}) can be rewritten as
\begin{equation} \label{eq:Xdisc}
X_h^n = S_{h,\Delta t}^n P_h x + \sqrt{\Delta t} \sum_{k=0}^{n-1} S_{h,\Delta t}^{n-k-1} 
(I + \theta \Delta A_h)^{-1} P_h \chi^{k+1}, \quad 0\leq n \leq N,
\end{equation}
\noindent where we have set for any $h>0$ and $N\geq 1$:
$$
S_{h,\Delta t} = (I + \theta \Delta t A_h)^{-1} (I - (1-\theta) \Delta t A_h ).
$$

{\bf Step 1:}
We introduce discrete and semi-discrete auxiliary schemes  which will 
be usefull for the proof of Theorem 
\ref{theo:weak}. 

First, for any $h>0$, let $\{X_h(t)\}_{t\in [0,T]}$ be 
the $V_h$-valued stochastic 
process solution of the following finite dimensional \spde{}
$$
\dd X_{h,t} + A_h X_{h,t} \, \dd t = P_h \, Q^{1/2} \dd W_t, \quad X_{h,0} = P_h x.
$$
It is straightforward to see that $X_{h,t}$ can be written as
\begin{equation} \label{eq:Xsem_sto_heat}
X_{h,t} = S_h(t) P_h x + \int_0^t S_h(t-s) P_h Q^{1/2} \dd W_s.
\end{equation}
The last stochastic integral is well defined since 
$t\mapsto \Tr{}((S_h(t) P_h Q^{1/2})^\star (S_h(t) P_h  Q^{1/2}))$ is integrable on $[0,T]$ .

We introduce also the 
following $V_h$-valued stochastic process
$$
Y_{h,t} = S_h(T-t) \, X_{h,t}, \quad t \in [0,T],
$$
\noindent which is solution of the following drift-free finite dimensional stochastic differential 
equation
\begin{equation}
\label{eq:Ysem_sto_heat}
\dd Y_{h,t} = S_{h}(T-t) P_h \, Q^{1/2} \dd W_t, \quad Y_{h,0} = S_h(T) P_h x.
\end{equation}
\noindent Its discrete counterpart is given by
\begin{eqnarray} \label{eq:Ydisc}
Y_h^n & = & S_{h,\Delta t}^{N-n} \, X_h^n, \quad 0\leq n \leq N,\\
      & = & S_{h,\Delta t}^N P_h x + \sqrt{\Delta t} \sum_{k=0}^{n-1} S_{h,\Delta t}^{N-k-1} 
(I + \theta \Delta t A_h)^{-1} P_h Q^{1/2} \chi^{k+1}.\nonumber
\end{eqnarray}

Eventually, we consider a time continuous interpolation of $Y_h^n$ which is the $V_h$-valued 
$\{{\cal F}_t\}_t$-adapted stochastic process $\widetilde{Y}_{h,t}$ defined by 
\begin{equation}\label{eq:Ytildesem}
\widetilde{Y}_{h,t} = S_{h,\Delta t}^N P_h x + 
\int_0^t \sum_{n=0}^{N-1} S_{h,\Delta t}^{N-n-1} (I + \theta \Delta t A_h)^{-1} {\mathbf 1}_n(s) 
P_h Q^{1/2} \dd W_s,
\end{equation}
\noindent where ${\mathbf 1}_n$ denotes the function ${\mathbf 1}_{[t_n,t_{n+1}[}$.

It is easy to see that for any $t\in [0,T]$ and $n$ be such that $t \in [t_n,t_{n+1}[$, we have 
$$
\widetilde{Y}_{h,t} = Y^n_h + S_{h,\Delta t}^{N-n-1}(I + \theta \Delta t A_h)^{-1} P_h Q^{1/2} (W_t - W_{t_n}).
$$

{\bf Step 2:} Splitting of the error.

Let now $\varphi\in C^2_b(H)$. The error $\E \, \varphi(X^N_h) - \E \, \varphi(X_T)$ can be 
splitted into two terms:
\begin{eqnarray} 
\label{eq:error}
\E \, \varphi(X^N_h) - \E \, \varphi(X_T) & = & \E \, \varphi(X^N_h) - \E \, \varphi(X_{h,T}) 
                                            + \E \, \varphi(X_{h,T}) - \E \, \varphi(X_T)  \\
                                          & = & A + B. \nonumber
\end{eqnarray}
\noindent The term $A$ contains the error due to the time discretization and will be estimated uniformly
with respect to $h$. The term $B$ contains the spatial error. 

%----------------------
% Erreur en temps
%----------------------

{\bf Step 3:} Estimate of the time discretization error. 

Let us now estimate the time error uniformly with respect to $h$. In order to do this, we consider the solution $v_h:V_h\rightarrow \R$ 
of the 
following deterministic finite dimensional Cauchy problem:
\begin{equation} \label{eq:kolmo_vh}
\left \{
\begin{array}{l}
\der{v_h}{t}{} = \frac 12 \Tr \left ((S_h(T-t) P_h Q^{1/2})^\star D^2 v_h 
                                     (S_h(T-t) P_h Q^{1/2})\right ), \\
v_h(0) = \varphi.
\end{array}
\right
.
\end{equation}
We have the following classical representation of the solution of (\ref{eq:kolmo_vh}) at any time $t\in[0,T]$ and for any $y \in V_h$:
\begin{equation} \label{eq:feynman}
v_h(T-t,y)   =   \E \, \varphi \left (  y + 
\int_t^T S_h(T-s) P_h \, Q^{1/2} \dd W_s \right ).
\end{equation}
It follows easily
\begin{equation} \label{eq:feynman_bis}
\|v_h(t)\|_{C^2_b(H)} \le   \|\varphi \|_{C^2_b(H)},\; t\in [0,T].
\end{equation}

Now, the estimate of the time error relies mainly on the comparison of 
It\^o formula applied successively to $t\mapsto v_h(T-t,Y_{h,t})$ and 
$t\mapsto v_h(T-t,\widetilde{Y}_{h,t})$.
First, by construction, $t\mapsto v_h(T-t,Y_{h,t})$ is a martingale. Indeed, It\^o formula 
gives 
$$
dv_h(T-t,Y_{h,t})= \left(Dv_h(T-t,Y_{h,t}),S_h(T-t)P_h Q^{1/2}dW_t\right).
$$
Therefore
$$
v_h(T-t,Y_{h,t})= v_h(T,S_h(T)P_h x) + \int_0^t\left(Dv_h(T-s,Y_{h,s}),S_h(T-s)P_h Q^{1/2}dW_s\right).
$$
Taking $t=T$ and the expectation implies
\begin{equation} \label{eq:feynman_1}
 \E \, \varphi(X_{h,T}) = v_h(T,S_h(T)P_h x).
\end{equation}
On the contrary,  $t\mapsto v_h(T-t,\widetilde{Y}_{h,t})$ is not a martingale. 
Nevertheless, applying It\^o formula gives, thanks to (\ref{eq:Ytildesem}),
{\setlength{\arraycolsep}{1mm}
%\begin{equation} \label{eq:feynman_2}
\begin{eqnarray} 
\E \, v_h(0,\widetilde{Y}_{h,T}) & = & \E \, v_h(T,\widetilde{Y}_{h,0}) 
- \E \, \int_0^T \der{v_h}{t}{}(T-t,\widetilde{Y}_{h,t}) \, \dd t \nonumber \\
 & + & \frac 12 \; \E \int_0^T \sum_{n=0}^{N-1} \Tr \bigg [ 
 \left ( S^{N-n-1}_{h,\Delta t} T_{h,\Delta t} P_h Q^{1/2}\right)^\star
 D^2v_h%(T-t,\widetilde{Y}_{h,t}) 
 \left ( S^{N-n-1}_{h,\Delta t} T_{h,\Delta t} P_h Q^{1/2}\right) 
 \bigg ] {\mathbf 1}_n(t)  \, \dd t,  \label{eq:feynman_2}
\end{eqnarray}
%\end{equation}
}
\hspace{-.15cm}where here and in equations (\ref{eq:feynman_3}), (\ref{eq:time_error}) below,  $D^2v_h$ is evaluated at $(T-t,\widetilde{Y}_{h,t})$.
Also we have set
$$
T_{h,\Delta t} = (I + \theta \Delta t A_h)^{-1}.
$$
Now, plugging (\ref{eq:kolmo_vh}) into (\ref{eq:feynman_2}) gives: 
{\setlength{\arraycolsep}{1mm}
%\begin{equation} \label{eq:feynman_3}
\begin{eqnarray}
\E \, \varphi(X_h^N)  & =  &  v_h(T,S_{h,\Delta t}^N P_h x) \nonumber \\
& + & \frac 12 \; \E \int_0^T \sum_{n=0}^{N-1} \Tr \bigg [ 
 \left ( S^{N-n-1}_{h,\Delta t} T_{h,\Delta t} P_h Q^{1/2}\right)^\star
 D^2v_h%(T-t,\widetilde{Y}_{h,t}) 
 \left ( S^{N-n-1}_{h,\Delta t} T_{h,\Delta t} P_h Q^{1/2}\right) 
\label{eq:feynman_3}\\
& - & 
 \left ( S_h(T-t) P_h Q^{1/2}\right)^\star
 D^2v_h%(T-t,\widetilde{Y}_{h,t}) 
 \left ( S_h(T-t) P_h Q^{1/2}\right) 
 \bigg ] {\mathbf 1}_n(t)  \, \dd t. \nonumber
\end{eqnarray}
%\end{equation}
}
\hspace{-.15cm}At last, the comparison between (\ref{eq:feynman_1}) and 
(\ref{eq:feynman_3}) leads to the following 
decomposition of the time error $A$
{\setlength\arraycolsep{2pt}
\begin{eqnarray} \nonumber
\E \, \varphi(X_h^N) & - & \E \, \varphi(X_{h,T})  =  v_h(T,S_{h,\Delta t}^N P_h x)  - v_h(T,S_h(T)P_h x) \\
 & + & \frac 12 \; \E \int_0^T \sum_{n=0}^{N-1} \Tr \bigg [ 
 \left ( S^{N-n-1}_{h,\Delta t} T_{h,\Delta t} P_h Q^{1/2}\right)^\star
 D^2v_h%(T-t,\widetilde{Y}_{h,t}) 
 \left ( S^{N-n-1}_{h,\Delta t} T_{h,\Delta t} P_h Q^{1/2}\right) 
\label{eq:time_error}\\
& - & 
 \left ( S_h(T-t) P_h Q^{1/2}\right)^\star
 D^2v_h%(T-t,\widetilde{Y}_{h,t}) 
 \left ( S_h(T-t) P_h Q^{1/2}\right) 
 \bigg ] {\mathbf 1}_n(t)  \, \dd t, \nonumber \\
 & =& I + II, \nonumber
\end{eqnarray}}
The term $I$ is the pure deterministic part of the time error. Thanks to the representation (\ref{eq:feynman}), we have
\begin{equation}\label{eq:I}
I \leq \| \varphi \|_{{\cal C}_b^1(H)} \| S_h(T) P_h - S_{h,\Delta t}^N P_h \|_{{\cal L}(H)} | x |.
\end{equation}
Thanks to 
(\ref{eq:implicit}) it is possible to bound $I$ uniformly with respect to $h$. More precisely, 
we have
\begin{eqnarray} \label{eq:time_error_det}
\|(S_h(N \Delta t) - S^N_{h,\Delta t})P_h\|_{{\cal L}(H)} & = & \sup_{i\geq 1} \left | e^{- N \lambda_{i,h}\Delta t } - F^N(\lambda_{i,h}\Delta t)  \right | \\
& \leq & \sup_{z\geq 0} \left | e^{- N z} - F^N(z)  \right | \nonumber \\
& \leq &  \frac{\kappa_4}{N} \leq \kappa_4 \Delta t,  \nonumber
\end{eqnarray}
\noindent for any $T\geq 1$ (e.g. see Theorem 1.1, p. 921 in \cite{leroux}). We have used the following notation:
$$
F(z) = \frac{1-(1-\theta) z}{1 + \theta z},\; z>0.
$$

Let us now see how to estimate the term $II$. First, using the symmetry of $D^2v_h$, 
we rewrite the trace term as
$$
\begin{array}{l}
\Tr{} \left ( \left ( S_{h,\Delta t}^{N-n-1} T_{h,\Delta t} P_h Q^{1/2} - S_h(T-t) P_h Q^{1/2} \right )^\star D^2 v_h \left ( S_{h,\Delta t}^{N-n-1} T_{h,\Delta t} P_h Q^{1/2} - S_h(T-t) P_h Q^{1/2} \right  ) \right ) \\
+ 2 \, \Tr{} \left ( \left ( S_h(T-t) P_h Q^{1/2} \right )^\star D^2v_h   \left ( S_{h,\Delta t}^{N-n-1} T_{h,\Delta t} P_h Q^{1/2} - S_h(T-t) P_h Q^{1/2} \right  ) \right )\\
 =  a_n(t) + b_n(t).
\end{array}
$$

Let now $\alpha>0$ and $\beta\geq 0$ such that (\ref{eq:trace}) and (\ref{eq:Q}) hold with $0<1-\alpha+\beta \leq 1$. Let 
$\gamma>0$ and $\gamma_1>0$ such that $0<\gamma<\gamma_1<1 - \alpha + \beta\leq 1$.

We first estimate the term $a_n(t)$. We use  (\ref{eq:trace_prop1}), (\ref{eq:trace_prop2}),
(\ref{eq:feynman_bis}) and (\ref{eq:HS})
to obtain
\begin{equation}\label{eq:a_time}
\begin{array}{ll}
a_n(t) &\le \|D^2v_h(T-t)\|_{C_b^2(H)}\\
& \times \Tr{} \left ( \left ( S_{h,\Delta t}^{N-n-1} T_{h,\Delta t} P_h Q^{1/2} - S_h(T-t) P_h Q^{1/2} \right )^\star \left ( S_{h,\Delta t}^{N-n-1} T_{h,\Delta t} P_h Q^{1/2} - S_h(T-t) P_h Q^{1/2} \right  ) \right )\\
& \le  \|\varphi\|_{C_b^2(H)} \| \left(S_{h,\Delta t}^{N-n-1} T_{h,\Delta t}  - S_h(T-t) \right) 
P_h Q^{1/2}\|_{{\cal L}_2(H)}^2\\
&\leq \|\varphi\|_{C_b^2(H)}  \| \left(S_{h,\Delta t}^{N-n-1} T_{h,\Delta t}  - S_h(T-t) \right) 
A_h ^{(1-\gamma_1)/2}P_h\|_{{\cal L}(V_h)}^2 \| A_h^{(\gamma_1-1)/2} P_h Q^{1/2}\|_{{\cal L}_2(H,V_h)}^2.
\end{array}
\end{equation}
Note that here $V_h$ is endowed with the norm of $H$. Let us set 
\begin{equation}
\label{Mn}
\begin{array}{ll}
M_n(t) & = \left \| \left(S_{h,\Delta t}^{N-n-1} T_{h,\Delta t}  - S_h(T-t) \right) 
A_h ^{(1-\gamma_1)/2}P_h \right \|_{{\cal L}(V_h)}\\
&=  \ds{\sup_{1\leq i\leq I(h)}  \left |\frac{F^{N-n-1}(\lambda_{i,h} \Delta t)}{1 + \theta \lambda_{i,h} \Delta t} - e^{-\lambda_{i,h}(T-t)} \right |\lambda_{i,h}^{(1-\gamma_1)/2}.}
\end{array}
\end{equation}
Using similar techniques as for the proof of the strong order (see e.g.  \cite{printems}), we have the following bound, for $n<N-1$
\begin{equation}\label{eq:Nn_1}
M_n(t) \leq C \Delta t^{\gamma/2} ((N-n-1)\Delta t)^{-(1-\gamma_1+\gamma)/2},
\end{equation}
\noindent where here and below $C$ denotes a constant which depends only on 
$\gamma_1 ,\; \gamma,\;\|\varphi\|_{C_b^2(H)},\; \|A^\beta Q\|_{{\cal L}(H)}$
and $ \Tr{}(A^{-\alpha})$. In particular these constants do not depend on $h$ or $\Delta t$.
The proof is postponed to the appendix at
the end of this article.

We then estimate the last factor in (\ref{eq:a_time}). Since $V_h\subset H$ and we have endowed 
$V_h$ with the norm of $H$, we may write 
$$
 \| A_h^{(\gamma_1-1)/2}P_h Q^{1/2}\|_{{\cal L}_2(H,V_h)}
 \le  \| A_h^{(\gamma_1-1)/2}P_h Q^{1/2}\|_{{\cal L}_2(H)}.
$$
Using (\ref{eq:HS}), we deduce
$$
 \| A_h^{(\gamma_1-1)/2} P_h Q^{1/2}\|_{{\cal L}_2(H,V_h)}
 \le  \| A_h^{(\gamma_1-1)/2} P_h A^{-\beta/2}\|_{{\cal L}_2(H)}
 \| A^{\beta/2} Q^{1/2}\|_{{\cal L}(H)}.
$$
Using  (\ref{interp_h}) and (\ref{eq:interp_T-2}), we have
$$
\begin{array}{ll}
\| A_h^{(\gamma_1-1)/2} P_h A^{-\beta/2}\|_{{\cal L}_2(H)}^2
&=\sum_{i\in \N} | A_h^{(\gamma_1-1)/2} P_h A^{-\beta/2} e_i|^2\\
&\le \sum_{i\in \N} | A_h^{-1/2} P_h A^{-\beta/2} e_i|^{2(1-\gamma_1)}|P_h A^{-\beta/2} e_i|^{2\gamma_1}\\
&\le \sum_{i\in•À?\N} |A^{-1/2-\beta/2} e_i|^{2(1-\gamma_1)}|A^{-\beta/2} e_i|^{2\gamma_1}\\
%&= \sum_{i\in\N} \lambda_i^{-(1+\beta)(1-\gamma_1)-\beta\gamma_1}\\
&=\sum_{i\in\N} \lambda_i^{-(1-\gamma_1+\beta)}.
\end{array}
$$
We deduce from $1-\gamma_1+\beta >\alpha$ and (\ref{eq:Q_lambda}) 
that
$$
\| A_h^{(\gamma_1-1)/2} P_h A^{-\beta/2}\|_{{\cal L}_2(H)}^2\le  
\lambda_1^{1-\gamma_1+\beta -\alpha}\| A^{\beta} Q\|_{{\cal L}(H)}^{1/2} \Tr (A^{-\alpha})
$$
Plugging this and (\ref{eq:Nn_1})  into (\ref{eq:a_time}) yields for $n<N-1$:
\begin{equation}\label{eq:a_time_II}
a_n(t) \leq C \Delta t^\gamma ((N-n-1)\Delta t)^{-(1-\gamma_1+\gamma)}.
\end{equation}
For $n=N-1$, we derive similarly, 
\begin{equation}\label{eq:a_N-1}
\begin{array}{ll}
a_{N-1}(t)&\le  \|\varphi\|_{C_b^2(H)}  \| \left( T_{h,\Delta t}  - S_h(T-t) \right) 
A_h ^{(1-\gamma_1)/2}\|_{{\cal L}(V_h)}^2 \| A_h ^{(1-\gamma_1)/2}P_h Q^{1/2}\|_{{\cal L}_2(H,V_h)}^2\\
&\le C \left(\| T_{h,\Delta t} A_h ^{(1-\gamma_1)/2}\|_{{\cal L}(V_h)}^2 +\|S_h(T-t) 
  A_h ^{(1-\gamma_1)/2}\|_{{\cal L}(V_h)}^2\right)\\
&\le C \left( \Delta t^{\gamma_1-1} + (T-t)^{\gamma_1-1}\right)\\
&\le  C(T-t)^{\gamma_1-1}.
\end{array}
\end{equation}
Concerning $b_n$, we write
$$
\begin{array}{ll}
b_n(t)&=2 \, \Tr{} \left ( \left ( S_h(T-t) A_h ^{(1-\gamma_1)/2} A_h ^{(\gamma_1-1)/2}  P_h
Q^{1/2} \right )^\star D^2v_h   \right.\\
& \hspace{4cm} \left. \left ( S_{h,\Delta t}^{N-n-1} T_{h,\Delta t}  - S_h(T-t)\right  ) A_h ^{(1-\gamma_1)/2} A_h ^{(\gamma_1-1)/2}  P_h Q^{1/2}  \right )\\
&\le \|\varphi\|_{C_b^2(H)} \|S_h(T-t) A_h ^{(1-\gamma_1)/2} \|_{{\cal L}(V_h)}
\| \left ( S_{h,\Delta t}^{N-n-1} T_{h,\Delta t}  - S_h(T-t)\right  ) A_h ^{(1-\gamma_1)/2} \|_{{\cal L}(V_h)}\\
&  \hspace{4cm} \times 
 \| A_h^{(\gamma_1-1)/2}P_h Q^{1/2}\|_{{\cal L}_2(H,V_h)}^2
\end{array}
$$
Using similar argument as above, we prove
$$
\|S_h(T-t) A_h ^{(1-\gamma_1)/2} \|_{{\cal L}(V_h)}\le C (T-t)^{(\gamma_1-1)/2}
$$
and, for $n<N-1$,
$$
\| \left ( S_{h,\Delta t}^{N-n-1} T_{h,\Delta t}  - S_h(T-t)\right  ) A_h ^{(1-\gamma_1)/2} \|_{{\cal L}(V_h)}
\le C \Delta t^{\gamma} ((N-n-1)\Delta t)^{-(1-\gamma_1+\gamma)}
$$
so that $n<N-1$:
\begin{equation}\label{eq:b_time_II}
b_n(t) \leq C \Delta t^\gamma ((N-n-1)\Delta t)^{-(1-\gamma_1+\gamma)}.
\end{equation}
For $n=N-1$, we have 
\begin{equation}\label{eq:b_N-1}
\begin{array}{ll}
b_{N-1}(t)&\le  \|\varphi\|_{C_b^2(H)}  \|S_h(T-t) A_h ^{(1-\gamma_1)/2} \|_{{\cal L}(V_h)} \\
&\times \| \left( T_{h,\Delta t}  - S_h(T-t) \right) 
A_h ^{(1-\gamma_1)/2}P_h\|_{{\cal L}(V_h)} \| A_h ^{(1-\gamma_1)/2}P_h Q^{1/2}\|_{{\cal L}_2(H,V_h)}^2\\
&\le C (T-t)^{\gamma_1-1}
\end{array}
\end{equation}
We are now ready to bound $II$ in (\ref{eq:time_error}). Indeed, (\ref{eq:a_time_II}),  
(\ref{eq:a_N-1}), (\ref{eq:b_time_II}), (\ref{eq:b_N-1}) imply
\begin{eqnarray}\label{eq:II}
II & \leq &C \int_{\Delta t}^T \sum_{n=1}^{N-2} \Delta t^\gamma ((N-n-1)\Delta t)^{-(1-\gamma_1+\gamma)} {\mathbf 1}_n(t) \, \dd t +\int_0^{\Delta t} C (T-t)^{\gamma_1-1} dt\\
   & \leq & C \Delta t^\gamma    \nonumber
\end{eqnarray}

Then, plugging (\ref{eq:I}) and (\ref{eq:II}) into (\ref{eq:time_error}) we obtain that 
\begin{eqnarray} \label{eq:time_end}
|A| &  \leq &  \kappa_4 \|D \varphi\|_{{\cal C}_b(H;H)} |x|\, \Delta t   +   \frac{C_3}{2}\, T^{\gamma_1 - \gamma} \Delta t^\gamma \\
     & \leq & C \, \Delta t^\gamma, \nonumber
\end{eqnarray}
\noindent for $T\geq 1$, $\Delta t\leq 1$.

%-----------------------
%  Erreur spatiale
%-----------------------

{\bf Step 4:}  Estimate of the space discretization error.

Let us now estimate the spatial error $B$. The method is essentially the same as above: we
use the Kolmogorov  equation associated to the transformed process $Y_t$. 

%For $n\geq 1$, 
%let us denote by $P_n$ the $H$-orthogonal projection onto $\mbox{Span}\{e_1,\dots,e_n\}$ 
%and by $V_n=P_n V$. 
We consider the following linear parabolic equation on $H$: 
\begin{equation} \label{eq:kolmo_vn}
\der{v}{t}{}(t,x) = \demi \Tr{}\Big ( D^2v(t,x) (S(T-t)Q^{1/2})(S(T-t)Q^{1/2})^\star \Big ), \quad t >0, \quad x \in H,
\end{equation}
\noindent together with the initial condition
$$
v(0,x) = \varphi(x), \quad x \in H,
$$
\noindent where $v$ is a real-valued function of $t$ and $x\in H$. We have the following representation of $v$ 
(see e.g. \cite{dpz:2ndorder}, chapter 3) at time $t\in[0,T]$ and any $y\in H$:
\begin{equation} \label{eq:rep_kolmo_vn}
v(T-t,y) = \E \, \varphi \left ( y + \int_t^T S(T-s) \, Q^{1/2} \dd W_s \right ).
\end{equation}
We apply the It\^o{} formula to $t\mapsto v(T-t,Y_t)$ and
$t\mapsto v(T-t,Y_{h,t})$. We substract the resulting equations and obtain
\begin{equation}\label{eq:space_error}
\left \{
\begin{array}{l} 
\E \varphi( X_T) - \E \varphi( X_{h,T})  =  v(T, S(T)x) - v(T,S_h(T)P_h x) \\
 \\
+\ds{\demi \, \E \int_0^T \Tr{}\bigg [ \Big ( S_h(T-t) P_h Q^{1/2} \Big )^\star D^2 v(T-t, Y_{h,t}) \Big (S_h(T-t) P_hQ^{1/2} \Big ) \bigg ] \, \dd t }\\
\\
-\ds{\demi \, \E \int_0^T \Tr{}\bigg [ \Big ( S(T-t) Q^{1/2} \Big )^\star D^2 v(T-t, Y_{h,t}) \Big (S(T-t) Q^{1/2} \Big ) \bigg ] \, \dd t .}
\end{array}
\right .
\end{equation}
The first term on the right hand side of (\ref{eq:space_error}) is the deterministic spatial 
error which can be bounded thanks to (\ref{eq:smooth}) (with $s=0$ and $q=2\gamma<2$) and (\ref{eq:rep_kolmo_vn}). We obtain:
\begin{equation} \label{eq:100}
\left | v(T, S(T)x) - v(T,S_h(T)P_h x) \right | \leq \kappa_1 \| \varphi\|_{C_b(H)^1} h^{2\gamma} T^{-\gamma} |x|.
\end{equation}
For the second term, we use  the symmetry of $D^2 v$ and write
$$
{\setlength\arraycolsep{2pt}
\begin{array}{ll}
& \Tr{} \bigg [ 

\Big ( S_h(T-t) P_h Q^{1/2} \Big )^\star D^2 v  \Big ( S_h(T-t) P_h Q^{1/2} \Big ) 
- \Big ( S(T-t) Q^{1/2}  \Big )^\star D^2 v \Big ( S(T-t) Q^{1/2} \Big ) \bigg ]\\
&\\
= & \Tr{} \bigg [

\Big ( S_h(T-t) P_h Q^{1/2} - S(T-t) Q^{1/2}  \Big )^\star D^2 v \Big ( S_h(T-t) P_h Q^{1/2} - S(T-t) Q^{1/2}  \Big )
\bigg ] \\

+ & 2 \, \Tr{} \bigg [

\Big ( S(T-t) Q^{1/2}  \Big)^\star D^2 v \Big ( S_h(T-t) P_h Q^{1/2} - S(T-t) Q^{1/2}  \Big ) \bigg ] \\
&\\
= & a + b,
\end{array}
}
$$
where here and below $D^2v$ is evaluated at $(T-t, Y_{h,t})$.

Let $\gamma>0$ be such that (\ref{eq:order}) holds and $\gamma_1>0$ such that $0<\gamma <\gamma_1< 1-\alpha +\beta\leq 1$.  Thanks to (\ref{eq:trace_prop2}) and (\ref{eq:rep_kolmo_vn}), we get the following bounds:
\begin{eqnarray*}
b & = & 2 \, \Tr{} \Big ( S(T-t) D^2 v (S_h(T-t)  P_h - S(T-t) )Q \Big ) \\
  & = & 2 \, \Tr{} \Big ( A^{\gamma_1-1-\beta} 
  A^{1-\gamma_1} S(T-t) D^2 v (S_h(T-t)  P_h - S(T-t) ) QA^\beta \Big ) \\
     & \leq &2\,  \| (S_h(T-t) P_h - S(T-t))\|_{{\cal L}(H)} \|\varphi\|_{{\cal C}_b^2(H)} 
     \|A^{1-\gamma_1} S(T-t)\|_{{\cal L}(H)} \| Q A^{\beta}\|_{{\cal L}(H)}\Tr (A^{1-\gamma_1-\beta}).
\end{eqnarray*}
Then, owing to (\ref{eq:regular0}), (\ref{eq:smooth}) (with $s=0$ and $q=2\gamma<2$), we obtain
\begin{equation}\label{eq:b_space}
b \leq C h^{2\gamma} (T-t)^{-(1+\gamma-\gamma_1)} 
\end{equation}
where again $C$ denotes a constant which depends only on 
$\gamma_1 ,\; \gamma,\;\|\varphi\|_{C_b^2(H)},\; \|A^\beta Q\|_{{\cal L}(H)}$
and $ \Tr{}(A^{-\alpha})$ but not on $h$ or $\Delta t$.

As regard $a$, we get first thanks to (\ref{eq:trace_prop2}) and (\ref{eq:rep_kolmo_vn}):
\begin{eqnarray}\nonumber
a & \leq & \| \varphi\|_{C_b^2(H)} \, \Tr{} \Big ( Q (S_h(T-t) P_h - S(T-t))^\star (S_h(T-t)P_h - S(T-t)) \Big )
\nonumber\\
&=& \| \varphi\|_{C_b^2(H)}\| Q^{1/2} (S_h(T-t) P_h - S(T-t))\|_{{\cal L}_2(H)}^2
\nonumber\\
   & \leq & \| \varphi\|_{C_b^2(H)}\| Q^{1/2} A^{\beta/2}\|_{{\cal L}(H)}^2
   \|A^{(1-\gamma_1)/2} (S_h(T-t) P_h - S(T-t))\|_{{\cal L}(H)}^2 \|A^{-(1-\gamma_1+\beta)/2}\|_{{\cal L}_2(H)}^2
   \nonumber\\
   &\le & C \|A^{(1-\gamma_1)/2} (S_h(T-t) P_h - S(T-t))\|_{{\cal L}(H)}^2
  ,\nonumber
 \end{eqnarray}
 \noindent where we have used (\ref{eq:Q_lambda}).
 If $\gamma_1+\gamma\ge 1$, we interpolate (\ref{eq:smooth}) with $q=(\gamma+\gamma_1-1)/\gamma$, $s=0$ and (\ref{eq:non-smooth0}) with 
 $s'=0$, $q'=1$
 and get
\begin{equation} \label{e3.63}
 \|A^{(1-\gamma_1)/2}  (S_h(T-t) P_h - S(T-t)) \|_{{\cal L}(H)} \leq C\, h^{\gamma} (T-t)^{-(1-\gamma_1+\gamma)/2}.
\end{equation}
If $\gamma_1+\gamma< 1$, we interpolate (\ref{eq:smooth}) with $q=0$, $s=0$ and (\ref{eq:non-smooth0}) with 
 $s'=0$, $q'=1$
 and get
\begin{equation} \label{e3.64}
 \|A^{(1-\gamma_1)/2}  (S_h(T-t) P_h - S(T-t)) \|_{{\cal L}(H)} \leq C\, h^{1-\gamma_1} (T-t)^{-(1-\gamma_1)}.
\end{equation}
We use again an interpolation argument to get
\begin{equation} \label{e3.65}
\begin{array}{l}
 \|A^{(1-\gamma_1)/2}  (S_h(T-t) P_h - S(T-t)) \|_{{\cal L}(H)}\\
  \leq C
 \| (S_h(T-t) P_h - S(T-t)) \|_{{\cal L}(H)}^{\gamma_1}\|A^{1/2}  (S_h(T-t) P_h - S(T-t)) \|_{{\cal L}(H)}^{1-\gamma_1}\\
 \le C \left(\|A^{1/2}  S_h(T-t) P_h \|_{{\cal L}(H)}+\|A^{1/2}S(T-t) \|_{{\cal L}(H)}\right)^{1-\gamma_1}\\
 \le C (T-t)^{(1-\gamma_1)/2},
 \end{array}
\end{equation}
thanks to (\ref{eq:regular0})  for $A$ and $A_h$ and (\ref{eq:interp_A}). A further interpolation 
between (\ref{e3.64}) and (\ref{e3.65}) gives
\begin{equation} \label{eq3.66}
 \|A^{(1-\gamma_1)/2}  (S_h(T-t) P_h - S(T-t)) \|_{{\cal L}(H)} \leq C\, h^{\lambda(1-\gamma_1)} 
 (T-t)^{-(1-\gamma_1)(\lambda+ (1-\lambda)/2)}.
\end{equation}
Taking $\lambda= \gamma/(1-\gamma_1)$ shows that (\ref{e3.63}) again holds for 
$\gamma_1+\gamma\le 1$.

We deduce 
\begin{equation} \label{eq:a_space}
a \leq C h^{2\gamma}  (T-t)^{-1+\gamma_1-\gamma}. \\
  \end{equation}

Plugging (\ref{eq:100}), (\ref{eq:b_space}) and (\ref{eq:a_space}) into (\ref{eq:space_error}) leads to, after time integration which is 
relevant since $1-\gamma_1+\gamma<1$:
\begin{equation}\label{eq:space_end}
|B| \leq C\; h^{2\gamma},
\end{equation}
\noindent for $T\geq 1$.

{\bf Conclusion:}
Gathering (\ref{eq:time_end}) and (\ref{eq:space_end}) in (\ref{eq:error}) ends the proof of Theorem
\ref{theo:weak}.

\section*{Appendix.}

We now prove the estimate (\ref{eq:Nn_1}) on $M_n$, $n<N-1$, defined in (\ref{Mn}).

We proceed as follows:
\begin{eqnarray}  \label{eq:Nn_0} 
M_n(t) & \leq & \sup_{1\leq i \leq I(h)} \left | \frac{F^{N-n-1}(\lambda_{i,h} \Delta t) - e^{-\lambda_{i,h}(T-t_{n+1})}}{1 + 
\theta \lambda_{i,h} \Delta t}\right | \lambda_{i,h}^{(1-\gamma_1)/2}\\
  &  +& \sup_{1\leq i\leq I(h)} \left | \frac{e^{-\lambda_{i,h}(T-t_{n+1})} - e^{-\lambda_{i,h}(T-t)}}{1+\theta \lambda_{i,h} \Delta t}\right |
  \lambda_{i,h}^{(1-\gamma_1)/2}\nonumber\\
  &  +& \sup_{1\leq i \leq I(h)} e^{-\lambda_{i,h}(T-t)} \left ( 1- \frac{1}{1 + \theta \Delta t \lambda_{i,h}}  \right ) \lambda_{i,h}^{(1-\gamma_1)/2}
   \nonumber \\
   & = & a_1 + a_2 + a_3. \nonumber
\end{eqnarray}
Thanks to (\ref{eq:time_error_det}) with $N$ replaced by $N-n-1$ we get 
\begin{eqnarray*}
a_1 & \leq & \frac{\kappa_4}{N-n-1} \sup_{i\geq 1} \left ( \frac{\lambda_{i,h}^{(1-\gamma_1)/2}}{1 + \theta \lambda_{i,h}\Delta t} \right )  \\
        & \leq & \frac{\kappa_4 \Delta t^{(\gamma_1 - 1)/2}}{(N-n-1)^{(1-\gamma_1+\gamma)/2}} 
             \sup_{i\geq 1} \left ( \frac{(\lambda_{i,h}\Delta t)^{(1-\gamma_1)/2}}{1 + \theta \lambda_{i,h}\Delta t} \right )\\
        & \leq & \frac{\kappa_4 \Delta t^{\gamma/2}}{((N-n-1)\Delta t)^{(1-\gamma_1+\gamma)/2}}.
\end{eqnarray*}
\noindent Indeed, since $(1-\gamma_1+ \gamma )/2< 1$, $(N-n-1) \geq (N-n-1)^{(1-\gamma_1 +\gamma)/2}$.
In the same way, we have 
\begin{eqnarray*}
a_2 &\leq &\sup_{i\geq 1} \left ( \frac{1 - e^{-(t_{n+1}-t) \lambda_{i,h}}}{(1+\theta \Delta t \lambda_{i,h})} 
\lambda_{i,h}^{(1-\gamma_1)/2} e^{-(N-n-1)\Delta t \lambda_{i,h}} 
\right ) \\
 & \leq & C(\gamma) \, \Delta t^{\gamma/2} \sup_{i\geq 1} \left (  
\lambda_{i,h}^{(1-\gamma_1 + \gamma)/2 } e^{-(N-n-1)\Delta t \lambda_{i,h}} 
\right ) \\
 & \leq & \frac{C(\gamma,\gamma_1) \, \Delta t^{\gamma/2}}{((N-n-1)\Delta t)^{(1-\gamma_1 + \gamma)/2}},
\end{eqnarray*}
\noindent where we have used that $|t-t_{n+1}|\le \Delta t$ and the inequality 
$|e^{-x} - e^{-y}|\leq C_\gamma |x-y|^{\gamma/2}$ and (\ref{eq:regular}).
Eventually, similar computations lead to
\begin{eqnarray*}
a_3 & \leq & \Delta t^{\gamma/2}\sup_{i\geq 1}
 \left (  \frac{\theta (\Delta t\lambda_{i,h})^{1-\gamma/2}}{1+\theta \Delta t \lambda_{i,h}} 
 \lambda_{i,h}^{(1-\gamma_1+\gamma)/2} e^{-(T-t)\lambda_{i,h}} \right )\\
    & \leq & C(\gamma,\gamma_1) \,\Delta t^{\gamma/2} (T-t)^{-(1-\gamma_1+\gamma)/2}  \\
    & \leq & C(\gamma,\gamma_1)\, \Delta t^\gamma ((N-n-1)\Delta t)^{-(1-\gamma_1+\gamma)/2},
\end{eqnarray*}
\noindent where we have used again the inequality (\ref{eq:regular}).

Gathering these three estimates in (\ref{eq:Nn_0}) yields (\ref{eq:Nn_1}), for $n<N-1$.

\end{document}